\documentclass{amsart}

\usepackage{graphicx}
\usepackage{pdfpages}
\usepackage{xcolor}
\usepackage{amsfonts}
\usepackage{amsmath}
\usepackage{amscd}
\usepackage{amssymb}
\usepackage{alltt}
\usepackage{url}
\usepackage{ellipsis}

\usepackage{isabelle}
\usepackage{isabellesym}
\usepackage{isabelletags}
\usepackage{comment}
\usepackage{pdfsetup}
\usepackage{railsetup}
\usepackage{wasysym}

\usepackage{tikz} % 
\usetikzlibrary{chains,shapes,arrows,%
 trees,matrix,positioning,decorations}

\newtheorem{theorem}{Theorem}[subsection]
\newtheorem{lemma}[theorem]{Lemma}

\newcommand{\ring}[1]{\mathbb{#1}}
\newcommand{\op}[1]{\hbox{#1}}

\title{The Formal Proof of the Kepler Conjecture: a critical retrospective}
\author{Thomas Hales}
\date{}
\date{Version: 9/22/2022}

\begin{document}

\maketitle

\begin{abstract} 
  The Kepler conjecture asserts that no packing of congruent balls in
  three-dimensional Euclidean space has density greater than that of the
  face-centered cubic packing.  In 1998, Sam Ferguson and I
  announced a computer-assisted proof of this conjecture.  Long delays
  in the refereeing process sparked a project to give a formal proof
  of the Kepler conjecture, which was completed in a large
  collaborative effort in 2014.  This article gives a critical
  reappraisal of that project.
\end{abstract}

\parskip=0.8\baselineskip
\baselineskip=1.05\baselineskip

\newenvironment{blockquote}{%
  \par%
  \medskip%
  \baselineskip=0.7\baselineskip%
  \leftskip=2em\rightskip=2em%
  \noindent\ignorespaces}{%
  \par\medskip}

\section{Introduction}

In the late sixteenth century, Sir Walter Raleigh asked his assistant
Thomas Harriot to compute the number of cannonballs piled in a pyramid
(in defense against the Spanish Armada).  Thomas Harriot obtained a
general formula as the binomial coefficient ${k+d-1\choose d}$ for the
number of congruent balls piled in a pyramid of any row size $k$ in
any dimension $d$.  In two dimensions $d=2$, Harriot's formula reduces
to the formula $k(k+1)/2$ for triangular numbers.  In three
dimensions, Harriot's formula reduces to a formula for pyramidal
numbers known in antiquity from a Sanskrit source.  The cannonball
pyramid is properly known as the face-centered cubic packing.

Harriot was an early proponent of the atomic theory, and he viewed
atoms as small balls.  Harriot and Johannes Kepler shared an interest in
optics, and they corresponded during the first decade of the
seventeenth century.  Harriot pushed the atomic theory in their
correspondence, but Kepler remained skeptical.

However, one evening crossing the Charles bridge in Prague, as it
started to snow, Kepler started to contemplate why snowflakes have six
sides.  His contemplations led to a booklet ``The Six-Cornered
Snowflake,'' an early essay in crystallography, which posits that
arrangements of minute particles explain observable symmetries in
nature.  In this booklet from 1611, Kepler described the face-centered
cubic packing and asserted that it will be the ``tightest possible, so
that in no other arrangement could more pellets be stuffed into the
same container.'' The notion of density was suggested to him by the
tightly packed seeds of a pomegranite, and the face-centered cubic
packing was suggested to him by the three-dimensional structure of
honeycomb cells.

David Hilbert made the Kepler conjecture part of the eighteenth problem, in
his influential list of problems proposed at the International
Congress in Paris in 1900.

The two-dimensional analogue of the Kepler conjecture asserts that the
densest packing of congruent circular disks in the plane is achieved
by the hexagonal lattice packing.  Although the first careful proofs
do not appear until much later, Axel Thue is often credited with a
nineteenth century proof.

The first detailed strategy to prove the Kepler conjecture was
formulated by L\'aszl\'o Fejes T\'oth in 1953 \cite{Fej53}.  He was also the
first to suggest a computer-assisted proof.

Over the years, there were some notorious false claims of a proof.
Buckminster Fuller, the creator of the geodesic dome, claimed a proof,
but his claim lacked any credibility.  Around 1990, a Berkeley
professor Wu-Yi Hsiang claimed a proof.  In rebuttal to his claimed
proof, in a heated debate, I published a paper containing explicit
counterexamples to his work \cite{Hales:1994:MI}.  For a few years,
the research area became toxic, until the debate was settled in my
favor.

Sam Ferguson and I annouced a solution to the Kepler conjecture in
August 1998.  Our solution appeared in a series of preprints posted to
the ArXiv \cite{hales1998kepler}.  Because of the conjecture's
notorious history, I had hoped that publication of our papers would be
swift. This was not to be.  A panel of twelve referees was assigned
the review. Review dragged on for years, until the referees became
exhausted and quit.  Because of lingering doubts, the full publication
of the proof did not occur until nearly eight years after submission.
Concurrent
with publication, the editors of the \emph{Annals of Mathematics}
issued a policy on computer-assisted proofs avowing that ``the
computer part may not be checked line-by-line, but will be examined
for the methods by which the authors have eliminated or minimized
possible sources of error\ldots'' \cite{annalsPolicy}.
% https://annals.math.princeton.edu/board

During those years in delay, growing out of my frustration, in an
effort to bypass the referees, I launched a project to give a formal
proof of the Kepler conjecture.  Although the project was born out of
frustration, in truth, it was also shaped by a broader vision of the
central importance of computer-assisted mathematics and formalization
in years to come.  In a large group collaboration, the formal proof of
the Kepler conjecture was completed in 2014.
This article offers some comments beyond the official account of the
formalization  \cite{hales2017formal}.

%The formal proof
%contains of about a half million lines of proof scripts.

\section{Statement}

The Kepler conjecture asserts that no packing of congruent balls in
Euclidean three-space has density greater than that of the
face-centered cubic packing.  This section gives the formal statement
of the theorem.

The statement is invariant under rescalings of the radius of the
balls, and it is convenient to normalize the radius of each ball at
$R=1$.  The congruent balls do not enter directly. Rather we identify
the packing with the set $V$ of centers of the balls.  The radius
normalization implies that every two distinct points in $V$ are
separated by distance at least the diameter $2$.  In the formal
statement, $V$ is a set of points in $\ring{R}^3$.

As was known to Carl Friedrich Gauss, the face-centered cubic packing
is the densest possible among lattices (that is the set $V$ of integer
linear combinations of a vector-space basis).

The face-centered cubic packing is not the only packing of congruent
balls that achieves the density $\pi/\sqrt{18} \approx 0.74$ of the
face-centered cubic packing.  The best known alternative is the
hexagonal close packing, which is also built by stacking hexagonal
planar layers of balls.  However, in the hexagonal close packing the
layers are shifted so that the balls of one layer nestle into a different
set of pockets of the previous layer.  Infinitely many different
close packings are obtained in this manner by shifting hexagonal
layers to nestle into different sets of pockets in the previous
layer.  The formal statement does not include a uniqueness statement.

The Kepler conjecture is defined formally in HOL Light by the statement:
\begin{verbatim}
`the_kepler_conjecture <=>
     (!V. packing V
            ==> (?c. !r. &1 <= r
                       ==> &(CARD(V INTER ball(vec 0,r))) <=
                           pi * r pow 3 / sqrt(&18) + c * r pow 2))`
\end{verbatim}
In HOL Light, ASCII replacements of logical symbols are used:
the character $!$ stands for a universal quantifier $\forall$,
and $?$ stands for an existential quantifier $\exists$.
Also, \& stands for the explicit embedding of the set $\ring{N}$ of natural
numbers into the field $\ring{R}$ of real numbers.  In a slightly
cleaned up version of the above, we might write
\begin{align*}
&\op{the\_kepler\_conjecture} \iff\\
&\quad (\forall\, V.\ \op{packing}\ V\\
&\quad\quad \Longrightarrow (\exists c.\ \forall r.\ 1 \le r \Longrightarrow
\op{card}(V \cap B(0,r)) \le
\frac{\pi r^3}{\sqrt{18}} + c r^2)).
\end{align*}
We took pains to remove all mention of measure from the theorem, which
might naturally appear when speaking of densities.

To simplify the main statement as much as possible, the Kepler
conjecture is formulated as a statement about a finite
packing of balls within a large spherical container of radius $r$.  By
considering finite packings in a bounded container, boundary effects
of the container come into play, and an error term $c/r$ enters into
the inequality (after dividing both sides of the inequality by $r^3$).  
The error $c/r$ tends to zero as the radius $r$ of
the container tends to infinity.
The condition $1\le r$ is an explicit way of encoding the
condition that the container's radius $r$ should be sufficiently large.  

To push simplicity further, even $\pi$ might have been eliminated from
the statement, if we had used a large cubical container rather than a
large spherical container.  (By stating the Kepler conjecture in terms
of the number of balls, rather than the volume of the balls, the
constant $\pi$ drops out.)

As stated in the theorem, the constant $c$ in the error term depends
on the packing $V$ and is not explicit.  As a matter of fact, the
constant $c$ can actually be chosen to be independent of $V$ and one
nonoptimal value has been calculated to be $c=24373$ in an informal
calculation \cite{Scharf}.
% Nadja Sharf.
% https://arxiv.org/pdf/1712.03568.pdf

\section{Proof Sketch}

\subsection{partitioning space}

The first major step of the proof reduces the theorem to a local
statement about a single ball and its immediately neighboring balls.
This reduction relies on a geometric partition of space into cells,
which is adapted to the packing $V$ (viewed as always as the set of
centers of balls).

There is an enormous amount of freedom in choosing the partition
of space.  Much of our effort in proving the Kepler conjecture
went into the design of a good partition of space.  

The most obvious choice of partition is the Voronoi cell
decomposition, which assigns (to each point in a packing $V$) the
polyhedral region of all points closer to that point than to any other
point in the packing $V$. This approach leads to a proof of the
two-dimensional analogue of the Kepler conjecture.  However, this
approach does not give the desired sharp bound in three dimensions,
because there are Voronoi cells (such as a regular dodecahedron) that
have smaller volume than the volume of the Voronoi cell (a rhombic
dodecahedron) of the face-centered cubic packing.

In 1953, L. Fejes T\'oth modified the approach slightly so that to
each point he associated a weighted linear combination of two layers
of truncated Voronoi cells.  In this form, he conjectured that the
corresponding bound on density is sharp.  Fejes T\'oth's conjecture is
plausible, but still unproved.  From a both combinatorial and
computational points of view, his conjecture is unwieldy.  A slight
modification of Fejes T\'oth's conjecture (that truncates Voronoi
cells by polyhedra rather than balls) can be formulated as a
first-order statement for real-closed fields.  (Truncation by balls
introduces volumes that depend on the transcendental $\pi$, but
polyhedral truncation gives volumes that can be encoded with algebraic
formulas.)  Thus, Tarski's quantifier elimination (QE) algorithm could
in principle be applied to decide the modified Fejes T\'oth
conjecture, which implies the Kepler conjecture (by a short argument).
However, QE is thoroughly impractical on a problem of this size.

Starting around 1989, I started toying with other geometric partitions
of space.  The Delaunay decomposition (which is dual to the Voronoi
cell decomposition) was one of my first attempts.  In 1994, I started
to experiment with hybrid partitions that combined the best features
of both the Voronoi and Delaunay partitions.  We called the cells of
the hybrid partitions \emph{decomposition stars}.  

The original proof from 1998 was based on a hybrid partition.  This
partition enabled computer calculations that ran in acceptable time.
However, as I revised the text in preparation for formalization, the text
started to grow to alarming lengths.  At one point, as I recall, the
revised proof text reached nearly 600 pages, and I feared that it would
reach 1000.  Much of this length stemmed from the complexity of the
decomposition star, when expressed more formally.

The situation abrubtly changed in 2007, when Christian Marchal
introduced a remarkable new geometric partition of space, adapted to
the Kepler conjecture.  His partition functioned in much the same way
as decomposition stars, but it had the effect of removing the
complexities from the text part of the proof (at the expense of more
cases and longer computer calculations).  At first, Marchal appeared
to claim to have a new proof of the Kepler conjecture, but inspection
of his work showed that the calculations were incomplete, and he
backed away from the claim \cite{marchal:2007}.

Nevertheless, his research contains ingenious ideas that brought a
beneficial transformation in the formalization project.  Starting in
2007, abandoning the decomposition stars, I adapted the proof and
calculations to Marchal's partition, and this resulted in a much
shorter final text.
% https://hal.archives-ouvertes.fr/file/index/docid/160350/filename/070615_Kepler.pdf

\subsection{outline}

An outline of the proof of the Kepler conjecture goes as follows.
First, as explained already, a geometric partition of space is used to
reduce the conjecture to a statement about a single congruent ball and
its immediate neighbors.  The combinatorial structure of the
neighboring congruent balls is encoded as a planar graph.  A list of
combinatorial properties is established that must hold for any planar
graph that encodes a counterexample to the Kepler conjecture.  A graph
with these properties is called a \emph{tame} planar graph.  A
computer search enumerates all tame planar graphs, up to isomorphism.

The remaining steps analyze one tame planar graph at a time, to rule
out each one as a potential counterexample.  The planar graph is a
combinatorial object.  To realize a counterexample, it would be
necessary to realize the graph as a geometric configuration in
Euclidean space, with nodes realized as centers of nonoverlapping
congruent balls and edges realized as metric conditions between
neighboring congruent balls.  This geometric configuration is highly
constrained by nonlinear constraints.  The constraints can be relaxed
to linear constraints, and linear programming is used to rule out each
potential counterexample.  When a single linear program does not
suffice, branch and bound algorithms are used.  When every potential
counterexample is eliminated, the Kepler conjecture is proved.

Throughout the proof, various nonlinear inequalities are used that are
checked by computer.  Thus, the full proof of the Kepler
conjecture consists of a 300-page text and three separate computer
programs: one to classify tame planar graphs, one to run linear
programs, and one to prove nonlinear inequalities.  The formalization
project verified both the text and the computer programs.

\subsection{revision}

Formalizers fall into two subgroups: the scribes and the revisionists.
The scribes see it their duty to make a faithful copy of the original
proof text.  They proudly set the original side-by-side against the
original to display their likeness. The revisionists, on the other
hand, strive to make a complete overhaul of the original proof, to
achieve something better.

I am a revisionist; the original proof of the Kepler conjecture was
completely transformed during formalization. Marchal cells replaced
decomposition stars, and we undertook a systematic investigation of
the properties of the cells.  Inspired by Georges Gonthier's
formalization of the Four-Color theorem, we replaced planar graphs
with hypermaps.  A new concept, called fans, was used to describe the
geometric realization of planar graphs.  The proof was modularized by
breaking it into largely independent sections.  It took an entire book
\emph{Dense Sphere Packings} to describe the transformed proof
\cite{DSP}.

Week by week, concerning my part in the collaboration, I asked which
would most advance the formalization: revising the mathematical proof
or formalizing the existing proof.  Until shortly before the end,
revision was always the better option.  One obvious reason for my
focus on mathematics is that I am a mathematician.  But another reason
is that our project (unlike others that formalize polished
second-generation proofs) was handling a freshly announced proof that
lacked the usual review from referees.

\section{Formalization}

I made the decision to use HOL Light, based on advice from Freek
Wiedijk.  HOL Light has its advantages: it has a well-developed theory
of real analysis; the system is easy to learn; and no boundaries
separate computer code from proof scripts.  A HOL Light proof script
is ordinary OCaml programming language code, with preprocessing of
expressions written between backquotes.  In particular, the tactic
language is a major functional programming language.  For example,
{\tt 0} denotes an integer in the programming language, whereas backquoted
{\tt `0`} denotes an (axiomatically constructed) natural number in HOL
Light.

The proof scripts amount to about a half million lines of code.  As
such, it stands as one of the largest formalization projects ever
completed.  The project took an estimated 20 work years to complete.
The main part of the verification can be processed on a laptop
computer in about five hours.  Three separate auxiliary verifications
are also required. They correspond to the three computer-assisted
portions of the proof: linear program verification, tame graph
classification, and nonlinear inequality verification.  The first two
auxiliary verifications take less than one day each.  The verification
of nonlinear inequalities takes about 5000 CPU hours in a cloud
computation.

We call the project \emph{Flyspeck}, for \emph{Formal Proof of the
  Kepler conjecture}.  All parts of the verification (including the
long nonlinear inequality verifications) have been checked by at least
three groups: the Flyspeck team before announcement, Josef Urban, and
an anonymous referee.  As part of a systematic audit of the
formalization, Mark Adams exported the main part of the verification
to the HOL Zero system and reconfirmed our results
\cite{adams2014flyspecking}.

The full Flyspeck project is freely available at our github page.%
\footnote{\url{https://github.com/flyspeck/flyspeck}}

\subsection{libraries}

John Harrison was the primary contributor to HOL Light libraries.  He
contributed libraries in trigonometry, multivariate differential
calculus, measure and integration in $\ring{R}^n$, convex sets, and
polyhedra \cite{DBLP:journals/jar/Harrison13}.  He contributed
powerful tactics such as the \emph{without-loss-of-generality} tactic,
which automates some of what mathematicians mean in proofs when they
write those words \cite{DBLP:conf/tphol/Harrison09}.

Aside from Harrison's libraries, most lemmas
were formulated minimally rather
than in generality, because our main interest was the narrow
goal of formalizing the Kepler conjecture rather than general
library development.  For example, the formalization avoids matrices,
integration on manifolds (such as a sphere), and the general Jordan
curve theorem (using only a special case).

After the project was completed, Harrison converted limits in his
analysis library from nets to filters to bring the library more in
line with other proof assistants.  This change required the proof
scripts to be updated.
Our project would benefit from further work to organize the project
and to make it available for other purposes (along the lines
of the Math Components Library that organized the Feit-Thompson
Odd Order theorem).

\subsection{international collaboration}

Formalization is an activity that naturally lends itself to
international collaboration.  The definitions and statements of lemmas
act as a specification of the project.  The top-level architecture of
the project shows how the main theorem follows from the lemmas.  Any
formal proof of a lemma -- no matter where, how, or by whom it is
produced -- advances the project.

The Flyspeck project was a large international collaboration.  The
main centers were the Isabelle group in Munich (directed by Tobias
Nipkow), the Hanoi group (managed by Adams and led by Hoang Le
Truong), my Pittsburgh group, and various independent contributors.
The collaboration was sustained by significant amounts of travel
between Munich, Pittsburgh, and Hanoi.  Some of the happiest news from
the Flyspeck project is that this collaboration worked.

Following its loose collaborative style, the project includes several
different tactic and proof styles: Harrison's efficient compact
proofs; Alexey Solovyev's port of SSReflect to HOL
Light;\footnote{SSReflect is a collection of tactics developed by
  Gonthier for the formal proof in Coq of the Four-Color theorem and
  used subsequently in the formalization of the Odd Order
  theorem. Solovyev ported SSReflect and several of its libraries to
  HOL Light \cite{solovyev2013thesis}.}; a large collection Hanoi
tactics (sometimes expressed in Vietnamese, such as the popular
\emph{nhanh} tactic, which is Vietnamese for \emph{fast}); and my own
idiosyncratic collection of tactic macros, which allowed me to quickly
generate long proof scripts with few keystrokes.

Occasionally loose collaboration became problematic, and the anarchy
was barely controlled.  For example, we had to grapple with three
libraries for lists that had subtle differences: an imported library from
Isabelle that came bundled with the tame graphs, an imported library
from Coq that arrived with Solovyev's port of SSReflect, and
Harrison's native HOL Light library.

\subsection{alignment between informal and formal}

The informal text (as given by Dense Sphere Packings - DSP) was
aligned with formal scripts through randomly generated seven letter
tracking codes.  As the informal and formal texts were edited and
rearranged, the codes gave persistent names to lemmas that enabled
tracking across different versions.  For example, the code OUIJTWY
appears in both the \LaTeX\ source file (DSP Lemma~2.25) and in the
HOL Light scripts for the following inverse trig identity (which
asserts that the two acute angles of a right triangle have sum
$\pi/2$).

\begin{lemma}{\tt [OUIJTWY]} If $y\in [-1,1]$, then
\[
\arccos(y) + \arctan_2\left(\sqrt{1-y^2},y\right) = \pi/2.
\]
\end{lemma}

\begin{verbatim}
let OUIJTWY = Trigonometry2.acs_atn2

let acs_atn2_t = 
  `!y. (-- &1 <= y /\ y <=  &1) 
       ==> (acs y = pi/(&2) - atn2(sqrt(&1 - y pow 2),y))`
\end{verbatim}

The \LaTeX\ source files also label each definition with the
corresponding term in the formal proof. For example, here is a snippet
of the \LaTeX\ source asserting the alignment
of the informal function $\arccos$ with the HOL Light function {\tt
  acs}:

\begin{verbatim}
  \begin{definition}[arccos]\guid{QZTBJMH}
  \formaldef{$\arccos$}{acs}
  \label{def:arccos}\formal{acs,\ ACS\_COS,\ COS\_ACS}
  ...
\end{verbatim}

This alignment is useful for documentation of the formal proof scripts
and for machine learning projects aimed at informal to formal
translation.

The \emph{DRY} acronym of software engineering stands for \emph{Don't
  Repeat Yourself}: each piece of knowledge should have a single
authoritative representation in software.  By creating both informal
and formal texts, we violate the principle, even if they are aligned.
Considerable effort was required to maintain alignment across edits.
Recently, my interests have shifted to controlled natural languages
(CNL), which are capable of automatically generating formal text from
informal \LaTeX\ source files.  In a second generation formal proof of
the Kepler conjecture, one of my requirements would be automated
alignment through a CNL.

\subsection{search}

Much of the time spent composing formal proof scripts is consumed by
search.  Each theorem has a name, and with few exceptions, every
invocation of the theorem requires the theorem to be cited by name.
Some systems have turned to hammers as a way to relieve the burden
that the cite-by-name paradigm imposes on users.

Roland Zumkeller and Harrison developed a heavily used search tool for
HOL Light.  Theorems can be searched by for name (according to any
regular expression) and by matching subterms (including type
information) in the theorem statement or conclusion.  The theorem name
includes the module name, permitting localized search within a given
module.  Boolean operations combine searches.  With practice, the tool
is very effective.

I recorded about ten thousand of my searches.  In many
cases, we can guess the lemma sought from the search term:

\begin{verbatim}
  `cat r []`              // concatenation with null list
  `?!`                    // all unique existence theorems
  `pi = x`                // all formulas for pi
  `a <= (b:num)`;`a <= b` // from natural number to real inequality
  `def "min"`             // find and print definition
  `name "WNW"`            // lookup by tracking code letters
  `a cross (b + c)`       // vector cross product linearity
\end{verbatim}

\section{Verifying the computer-assisted proof}

\subsection{tame graphs}

The first major success in the formalization of the Kepler conjecture
was the Gertrud Bauer-Nipkow formalization of the classification of
tame planar graphs in Isabelle \cite{DBLP:conf/cade/NipkowBS06}.  (As
mentioned above, the word \emph{tame} has a technical meaning in the
proof.  A tame graph is one that satisfies a long list of technical
conditions that means that the graph encodes the combinatorial
structure of a potential counterexample to the Kepler
conjecture. Thus, the classification of tame planar graphs gives the
classification of possible counterexamples.)

The Bauer-Nipkow formalization was crucial in launching the Flyspeck
project and in galvinizing interest.  However, as the earliest part of
the project, it is currently the least compatible piece of the entire
Flyspeck project.  A reimplementation of their work is sorely needed.

The Bauer-Nipkow classification is the only part of the formalization
done in Isabelle.  There have been many discussions about merging all
parts of the Flyspeck project into the same proof assistant, but this
is still work in progress.  Currently the statement of the
classification theorem is translated by hand from Isabelle to HOL
Light, and the translated statement is accepted as an unproved
postulate in the HOL Light part of the formalization.  The Isabelle
proof is not translated into HOL Light.  This hand translation is the
weakest link in the entire Flyspeck project.

As a step towards verifying the classification in HOL Light, in
unpublished work from 2018, Wiedijk has produced an automated
translation of the tame graph classification from Isabelle to OCaml,
which Harrison then (semi-manually) translated to definitions in HOL Light.
(Someone needs to check the compatibility of Harrison's translation with mine.)
Solovyev then created a tool that translates HOL Light equational theorems
into optimized executable OCaml code that computes using primitive inference rules.%
\footnote{\url{https://github.com/monadius/compute-hol-light}} His
experiments suggest that a verification of the tame graph
classification might be done in HOL Light in a 3000 hours computation.
(Isabelle is much faster because it exports code to ML, then executes
without the overhead of HOL primitive inferences.)
% Private communication 2018, R--Flyspeck.

A major inefficiency in the tame graph classification comes through
the definition of \emph{planar graph}.  The definition of planar graph
is treated as a black box in the Isabelle formalization: a planar
graph is any graph generated by an algorithm we will call $X$, where
algorithm $X$ is an edge-adding algorithm that I designed and
implemented in Mathematica (and later in Java) to generate all planar
graphs (with constraints on the number of faces, nodes, and degrees of
the nodes) around 1993.  The Bauer-Nipkow procedure filters the list
of planar graphs (generated by the black box algorithm) to produce the
sublist of tame graphs.

One of the cardinal rules of code verification is never
to verify code that has not been written with verification in mind.
The cardinal rule is violated for algorithm $X$.  It was a very
unpleasant experience for me to verify correctness in HOL Light of a
HOL Light translation of an Isabelle translation of an ML translation
of an old Java program written without formalization in mind.  In
retrospect, there would have been much better ways to proceed.  In
particular, the data structures used to represent planar graphs in
algorithm $X$ are entirely inappropriate for formalization.  (The data
encoding a graph is highly redundant, and the formalization has to
give a series of lemmas proving that the redundant data remains
in sync.)  Moreover, other planar graph generating algorithms
are available, such as \emph{plaintri}.

\subsection{linear programs}

Steven Obua's thesis initiated the verification of linear program for
the Flyspeck project
\cite{DBLP:conf/tphol/Obua05}\cite{DBLP:journals/amai/ObuaN09}.  The
verification of linear programs was completed by Solovyev's thesis
\cite{solovyev2013formal}.  Solovyev has developed a stand-alone tool
for the verification of linear programming certificates in HOL Light.
His tool is highly optimized \cite{solovyev2011efficient}.

The linear programming part of the proof is not simply a matter of
running linear programs.  The linear programs are generated by a
complex procedure from the combinatorial structure of the tame planar
graphs.  The output from the linear program must be translated back
into bounds on an underlying nonlinear optimization generated from the
graph.

\subsection{nonlinear inequalities}
Although others made contributions, Solovyev deserves enormous credit
for the verification of the computer-assisted portions of the proof
(the linear programs and nonlinear inequalities)
\cite{solovyev2013formal}.

It was understood from the very beginning of the formalization of the
Kepler conjecture that the largest challenge would be the
formalization of the nonlinear inequalities.  Specifically,
floating-point arithmetic is slower by orders of magnitude when
executed as logical rules in HOL Light than in a native processor.

The proof of the Kepler conjecture relies on about a thousand
nonlinear inequalities over the real numbers.  The inequalities
have the general form 
\begin{equation}\label{eqn:ineq}
\forall x \in D,\quad f_1(x) < 0 \lor f_2(x) < 0 \lor \cdots \lor f_k(x) < 0.
\end{equation}
Each function $f_i : R_i \to \ring{R}$ is defined on a subset $R_i$ of the
domain $D\subset \ring{R}^n$.
The inequality (\ref{eqn:ineq}) means more precisely that 
at each point $x\in D$, there exists $1\le i\le k$ such that
$x \in R_i$ and $f_i(x) < 0$.  The functions are generally analytic.

The box-shaped domain $D$ is a product of compact real intervals
$[a,b]$, and $n$ is small (usually $n\le 6$).  The dimension $n=6$ of
the domain $D$ occurs naturally, because a simplex in three dimensions
is determined by its six edges.  These inequalities are proved by
computer using interval arithmetic.

There was an extraordinary amount of freedom in the creation
of a finite set of nonlinear inequalities that collectively imply
the Kepler conjecture.  In no sense is the current collection
optimal.  Much further research is justified.

Starting in January 1994, I directed all my efforts toward a proof of
the Kepler conjecture.  My decision to put full energies behind the
conjecture was shaped by a key factor among others: I first learned in
1993 that interval arithmetic provided a reliable way to prove
nonlinear inequalities.  Ferguson joined me, and much of our research
between 1995-1998 went into the development of efficient algorithms
for proving nonlinear inequalities.  Eventually, we settled on using
automatic differentiation to compute second order Taylor expansions of
the nonlinear functions with certified bounds on the error terms.
When a Taylor expansion fails to give sufficient accuracy on a given
box-shaped domain, we branch and bound (by adaptively bisecting the
box and recursing on the smaller domains).  When a partial derivative
has fixed sign over the box, the function is monotone, and the
inequality is implied by a restricted inequality over a face of the
box, and the dimension drops.  (We do not claim originality in these
techniques, but considerable work went into efficient implementation.)

Ferguson and I gave independent implementations of the algorithms (in
C and C++) that we used to cross-check our code for bugs. Ferguson
finished his thesis in 1997, but returned to University of Michigan
for an extended visit in 1998 to prove the final outstanding nonlinear
inequalities as we prepared to announce the solution.

To prepare for formalization, I ported the code to OCaml, and spent
months reworking the nonlinear inequalities until the programs
terminated in about 10 hours (of intensive floating-point
calculations).  In the end, the collection of nonlinear inequalites
used for the formalization is completely different from the collection
of inequalites that was used in the original 1998 proof.  The
algorithms remained very close to those used earlier.

Solovyev was entirely responsible for formalization in HOL Light.  An
enormous amount of code optimization was necessary to bring the formal
calculation within the reach of a large cloud computation.  Some of
Solovyev's optimizations included performing arithmetic in a large
base rather than base $2$, building massive addition and
multiplication tables inside HOL-Light so that each arithmetic
operation could be done by a simple look-up, precomputing the number
of digits precision needed for each calculation and using the minimum
possible, computing all the decision points informally and then
replaying the decision-point scripts in the formal verification, and
caching intermediate expressions to avoid recomputation.  He carefully
crafted certain inner loops to use the minimal number of primitive HOL
Light operations.

Solovyev developed and documented an independent tool for automated
verification of nonlinear inequalities that can be used separately
from the Flyspeck project.  His tool has excellent performance
compared with tools in other systems.

Further major optimizations are certainly possible, and it would be an
interesting future project, to optimize to the point of moving the
entire formalization of the Kepler conjecture from the cloud to a
single laptop computer.  Some suggestions for further optimization
include implementing backchaining algorithms to compute partial
derivatives, exploiting common structures shared among inequalities to
verify inequalities in batches, precomputing second derivative bounds
on functions that appear in multiple inequalities (this was an
important technique in the 1998 proof that we later abandoned for no
good reason; in current the implementation, second derivative bounds
are wastefully recomputed every time the algorithm branches into two
sub-boxes), redesigning the collection of nonlinear inequalities,
precomputing lookup tables for inverse trig and other transcendental
functions, automatic generation of linear combinations (if $\sum_i a_i
f_i(x) < 0$ with $a_i\ge 0$, then for some $i$, we have the desired
conclusion $f_i(x) < 0$ in Equation (\ref{eqn:ineq})), better
decisions about partitioning boxes into sub-boxes, and generally any
number of techniques for global optimization developed during the past
twenty years (since our proof still relies on 1990s era technology).

The lion's share of the time is spent on computations near the
inequality \emph{``hot spots''}: small boxes where the functions
$f_i(x)$ are nearly $0$.  Speeding up the hot-spot computations would
bring significant gains.  One technique that should help would be to
center the Taylor expansion at the spot where $f_i(x)$ is maximized,
rather than at the center of the box (as we do now).

Machine learning might be tried.  We can view the space of proofs of
the Kepler conjecture as parametrized by various finite collections
$C$ of nonlinear inequalities.  To each collection $C$ there is a cost
as measured by the time required to verify the collection of
inequalities.  We wish to learn a collection $C$ that minimizes cost.

\newpage

%\printbibliography

\bibliography{refs} 
\bibliographystyle{alpha}

\end{document}